\documentclass[12pt]{article}
\usepackage[utf8]{inputenc}
\usepackage[T1]{fontenc}
\usepackage{geometry}
\usepackage{amsmath, amssymb, amsthm, amsfonts}
\usepackage{hyperref}
\hypersetup{pdftitle={Approximate Computation via Decision-Theoretic Simulability},pdfauthor={Deniz Akdemir}}
\usepackage{graphicx}
\usepackage{booktabs}
\usepackage{natbib}
\usepackage{algorithm}
\usepackage{algorithmic}
\usepackage{xcolor}

\newtheorem{theorem}{Theorem}[section]

\newtheorem{proposition}{Proposition}[section]

\newtheorem{definition}{Definition}[section]

\title{Approximate Computation via Decision-Theoretic Simulability}
\author{Deniz Akdemir}
\date{November 16, 2026}

\begin{document}

\maketitle

\begin{abstract}
Approximate algorithms are often evaluated through numerical error or performance on a selected task. We study a complementary question: when can an efficient simulator reproduce a target statistical experiment well enough to preserve every bounded downstream decision? We formulate simulation error using Le Cam deficiency restricted to computationally admissible kernels, and distinguish a finite-size error envelope from a guarantee witnessed by one uniform algorithm across input sizes. We give explicit risk-transfer and composition bounds and explain the direction of simulation required for approximate computation. A representation trade-off follows from the classical triangle inequality for deficiency: source reconstruction error plus per-parameter representation mismatch bounds the target-to-source deficiency from below. An exactly solvable binary example illustrates the trade-off. We also separate task-specific diagnostics and optimization surrogates from valid deficiency certificates. The results organize established experiment-comparison tools for computational approximation; they do not establish a new separation of complexity classes or a general learning-consistency theorem.

\end{abstract}

\section{Introduction}
An approximate computation may change a numerical output substantially while preserving the action taken by a particular decision rule. Conversely, small error on one task need not preserve information needed by another. These observations motivate evaluating approximate computation through statistical experiment comparison: does the available output permit simulation of the experiment required by a downstream user?

Le Cam's theory provides a language for this question \citep{lecam1964sufficiency,lecam1986asymptotic,torgersen1991comparison}. Its operational object is directional deficiency, defined by the best parameter-independent Markov kernel transforming one family of distributions into another. Restricting the admissible kernels introduces an explicit computational constraint. A useful certificate consists of an implementable kernel and a uniform error bound, rather than an assertion that such a kernel should exist.

This perspective complements approximation algorithms, randomized computation, and task-specific evaluation. It does not require treating classical computation or information theory as exclusively concerned with exactness. Nor does performance of a heuristic on selected instances establish a worst-case simulation guarantee for an NP-hard problem.

\paragraph{Scope and contributions.}
Section~\ref{sec:foundations} defines the statistical and computational objects and states a risk-transfer bound, including its runtime and approximation qualifications. Section~\ref{sec:representation} explains the direction required for representation recovery. Section~\ref{sec:computation} specifies computational experiments and proves two distinct composition bounds. Section~\ref{sec:transfer} derives an invariance--reconstruction trade-off and an exact finite example. Section~\ref{sec:optimization} distinguishes certificates from optimization surrogates. Section~\ref{sec:discussion} discusses limitations.

The risk-transfer and triangle-inequality arguments are consequences of classical experiment comparison theory. The contribution here is their explicit organization around computational approximation, including the uniformity conditions and failure modes of proposed certificates. No claim of priority for computationally restricted experiment comparison is made.

\section{Experiments and Computationally Admissible Simulation}\label{sec:foundations}
We work with finite observation and action spaces unless stated otherwise. This avoids measurable-selection qualifications; the arguments extend whenever the relevant kernels and compositions are well defined. An experiment is $\mathcal E=(\mathcal X,\Theta,\{P_\theta\})$. A second experiment $\mathcal F=(\mathcal Y,\Theta,\{Q_\theta\})$ uses the same parameter space. A kernel $K:\mathcal X\rightsquigarrow\mathcal Y$ cannot depend on the unknown $\theta$.

Use the convention
\[
 \|P-Q\|_{\rm TV}=\sup_A|P(A)-Q(A)|=\tfrac12\sum_x|P(x)-Q(x)|.
\]
For $f$ taking values in $[0,1]$, expectation differences are at most this distance. For signed $f$ with $|f|\le1$, the general bound has an additional factor of two.

\begin{definition}[Directional deficiency]
\[
 d(K;\mathcal E,\mathcal F)=\sup_\theta\|P_\theta K-Q_\theta\|_{\rm TV},\qquad
 \delta(\mathcal E,\mathcal F)=\inf_K d(K;\mathcal E,\mathcal F).
\]
Small $\delta(\mathcal E,\mathcal F)$ means that observations from $\mathcal E$ can simulate $\mathcal F$. The symmetric distance is $\Delta(\mathcal E,\mathcal F)=\max\{\delta(\mathcal E,\mathcal F),\delta(\mathcal F,\mathcal E)\}$.
\end{definition}

\subsection{Computational restrictions and uniformity}
For computational statements, specify sequences $(\mathcal E_n,\mathcal F_n)$, their encodings, and polynomial bounds on observation lengths. A uniform simulator is one probabilistic Turing machine that receives $(1^n,x)$ and produces an encoded target observation using random bits. Fix its time budget, for example a polynomial $p(n)$, before optimization. Write $\mathcal K_p$ for the class of machines satisfying that budget on every valid input. Define the size-$n$ envelope
\[
 \delta_{p,n}(\mathcal E_n,\mathcal F_n)=\inf_{K\in\mathcal K_p}
 \sup_{\theta\in\Theta_n}\|P_{n,\theta}K_n-Q_{n,\theta}\|_{\rm TV}.
\]
Then $\delta(\mathcal E_n,\mathcal F_n)\le\delta_{p,n}(\mathcal E_n,\mathcal F_n)$. More permissive budgets cannot increase this envelope.

The envelope is a pointwise infimum. A different machine may approach it at each $n$; it does not by itself supply one machine with that performance for all sizes. In particular, taking an unrestricted infimum over all polynomial-time machines separately at every finite size can admit machines with size-specific hardcoded tables and does not establish computational hardness. We therefore state efficient guarantees using an explicit uniform machine $K$ and its whole error sequence $d_n(K)$, with its runtime bound. For a tolerance input $\varepsilon$, any polynomial dependence on $1/\varepsilon$ must also be stated. We use computational deficiency as this constrained simulation criterion, rather than treating the notation $\delta_{\rm poly}$ alone as a complexity separation.

\subsection{Risk transfer}
A randomized rule $\rho:\mathcal Y\rightsquigarrow\mathcal A$ has risk
$R_\mathcal F(\theta,\rho)=\mathbb E_{Q_\theta\rho}\ell(\theta,A)$ for loss $\ell:\Theta\times\mathcal A\to[0,1]$.

\begin{theorem}[Witnessed risk transfer]\label{thm:risk}
If a kernel $K$ satisfies $d(K;\mathcal E,\mathcal F)\le\varepsilon$, then the composed rule $K\rho$ satisfies
\[
 \sup_\theta|R_\mathcal E(\theta,K\rho)-R_\mathcal F(\theta,\rho)|\le\varepsilon.
\]
If $K$ and $\rho$ are uniform polynomial-time procedures with polynomially bounded intermediate lengths, so is $K\rho$.
\end{theorem}
\begin{proof}
For fixed $\theta$, put $f_\theta(y)=\mathbb E_{\rho(y)}\ell(\theta,A)\in[0,1]$. The expectation bound for TV gives the inequality. Runtime is the sum of the two runtimes after substituting the intermediate-length bound. Efficient simulation alone does not make an inefficient target rule efficient.
\end{proof}
For every $\eta>0$, an infimum over a nonempty kernel class supplies a kernel with error at most its infimum plus $\eta$. Without attainment, the corresponding conclusion is $\delta+\eta$, not existence of a rule attaining an exact $\delta$ bound. For losses in $[a,b]$, multiply the error bound by $b-a$.

\section{Representations and the Direction of Recovery}\label{sec:representation}
Let $T:\mathcal X\rightsquigarrow\mathcal Z$ be a representation and let $\mathcal E T=\{P_\theta T\}$. Always $\delta(\mathcal E,\mathcal E T)=0$: applying $T$ simulates its own output. This forward direction says nothing about information preservation. Recovery requires small $\delta(\mathcal E T,\mathcal E)$.

\begin{proposition}[Recovery certificate]
If $D:\mathcal Z\rightsquigarrow\mathcal X$ satisfies $\sup_\theta\|P_\theta TD-P_\theta\|_{\rm TV}\le\varepsilon$, then every bounded-loss rule on the original experiment can be implemented after representation with risk difference at most $\varepsilon$.
\end{proposition}
\begin{proof}
Apply Theorem~\ref{thm:risk} to source $\mathcal E T$, target $\mathcal E$, and kernel $D$.
\end{proof}
A sufficient statistic admits exact recovery under the usual finite-space factorization conditions. Approximate recovery is weaker and needs an explicit bound. We do not assert a strict hierarchy from likelihood distortion through pairwise testing to global deficiency. Forward pairwise deficiency after applying a statistic is already zero. Even meaningful reverse pairwise bounds can use different recovery kernels for different pairs; a global certificate requires one kernel valid for all parameters. Also, writing $T^{-1}(z)$ inside a likelihood ratio is not a definition for a many-to-one representation.

\subsection{One successful task is insufficient}
Consider a parameter $\theta=(u,v)\in\{0,1\}^2$, a full observation $X=(u,v)$, and representation $Z=u$. Prediction of $u$ has zero error before and after compression, but $v$ cannot be recovered uniformly. Any decoder, at fixed $u$, must approximate two distinct point masses corresponding to $v=0,1$. Their TV distance is one, so the triangle inequality forces worst-case reconstruction error at least $1/2$. Uniformly guessing $v$ attains $1/2$. Thus
\[
 \delta(\mathcal E T,\mathcal E)=1/2
\]
despite perfect performance on the selected task. Empirical agreement on a risk stratification rule is therefore not an upper bound on full deficiency. A task-restricted evaluation may still be useful, but its decision class, sampling uncertainty, and distinction from uniform experiment recovery must be explicit.

\section{Computational Experiments and Composition}\label{sec:computation}
\subsection{An instance experiment and its target}
For a deterministic function $f_n$ on encoded instances, define the available experiment $\mathcal I_n=\{\delta_x:x\in\Theta_n\}$ and ideal target $\mathcal F_n=\{\delta_{f_n(x)}:x\in\Theta_n\}$. An algorithm $A$ receives $x$, not the answer. Its simulation error is
\[
 d_n(A)=\sup_x\|\mathcal L(A(1^n,x))-\delta_{f_n(x)}\|_{\rm TV}
       =\sup_x\Pr\{A(1^n,x)\ne f_n(x)\}.
\]
Consequently, full-TV simulation of a point-mass target requires high-probability exact output. Near-optimal objective value or small numerical error alone is insufficient. Approximate decision validity requires a specified coarser target or a specified restricted loss class; it cannot be deduced merely by changing the word ``exact'' to ``semantic.'' An arbitrary unconstrained statistical kernel can compute $f_n$ by lookup, which is why runtime restrictions matter.

\subsection{A qualified LeCam-P convention}
For a specified pair of experiment sequences and fixed tolerance schedule $\varepsilon_n$, call the pair efficiently simulable if there exist one uniform randomized algorithm $K$ and one polynomial $p$ such that its runtime is at most $p(n)$ and $d_n(K)\le\varepsilon_n$ for every $n$. This is the precise pair-relative meaning we attach to the earlier name LeCam-P. It is not a problem class independent of the observation model, target, and tolerance. We establish no containment or separation involving $P$, $NP$, or randomized complexity classes. In particular, membership for TSP or protein-structure prediction is not established here.

\subsection{Chains of experiments}
\begin{proposition}[Directional composition]\label{prop:chain}
For experiments on a common parameter space,
\[
 \delta(\mathcal E_0,\mathcal E_k)\le\sum_{i=1}^k\delta(\mathcal E_{i-1},\mathcal E_i).
\]
If explicit kernels $K_i$ have errors at most $\varepsilon_i$, their composition has error at most $\sum_i\varepsilon_i$. For fixed $k$, uniform polynomial-time kernels with polynomial intermediate-length bounds compose in polynomial time.
\end{proposition}
\begin{proof}
For two kernels, TV contraction and the triangle inequality give
\[
 \|P_\theta K_1K_2-Q_{2,\theta}\|_{\rm TV}
 \le\|P_\theta K_1-Q_{1,\theta}\|_{\rm TV}
    +\|Q_{1,\theta}K_2-Q_{2,\theta}\|_{\rm TV}.
\]
Take suprema and iterate. For the infimum statement, choose each kernel within $\eta/k$ of its deficiency and let $\eta\downarrow0$. For a growing number of stages, a polynomial bound on the total runtime and all intermediate lengths must be verified separately.
\end{proof}

\subsection{Adaptive pipelines need conditional bounds}
\begin{proposition}[Transcript perturbation]\label{prop:adaptive}
Suppose ideal and approximate pipelines share an initial law and conditional response kernels $M_{i,\theta}(\cdot\mid h)$ and $\widetilde M_{i,\theta}(\cdot\mid h)$ on the same response spaces. If
\[
 \sup_{\theta,h}\|M_{i,\theta}(\cdot\mid h)-\widetilde M_{i,\theta}(\cdot\mid h)\|_{\rm TV}\le\varepsilon_i,
\]
then the TV distance between their full transcript laws is at most $\min\{1,\sum_i\varepsilon_i\}$. The same bound holds after a common output-processing kernel.
\end{proposition}
\begin{proof}
Write $H_i,\widetilde H_i$ for the prefix laws. Adding and subtracting $\widetilde H_{i-1}M_i$ gives
\[
 \|H_i-\widetilde H_i\|_{\rm TV}
 \le\|H_{i-1}-\widetilde H_{i-1}\|_{\rm TV}+\varepsilon_i.
\]
Induction proves the sum bound; TV is at most one. Common postprocessing contracts TV.
\end{proof}
Unconditional marginal deficiencies alone do not establish this adaptive guarantee. Conditional kernels must be specified even on histories reached only by the approximate process. If response laws are first simulated, the bounds above must apply to the actual simulated conditional responses.

\section{Invariance and Reconstruction}\label{sec:transfer}
Let source and target experiments $\mathcal E_S,\mathcal E_T$ share $\Theta$, and let representation kernels produce observations in a common space $\mathcal Z$. Write $\mathcal F_S=\{Q^S_\theta\}$ and $\mathcal F_T=\{Q^T_\theta\}$ for the represented experiments, and
\[
 e=\sup_\theta\|Q^T_\theta-Q^S_\theta\|_{\rm TV}.
\]
\begin{theorem}[Invariance--reconstruction trade-off]\label{thm:transfer}
\[
 \delta(\mathcal F_S,\mathcal E_S)+e\ge\delta(\mathcal E_T,\mathcal E_S).
\]
\end{theorem}
\begin{proof}
The target representation gives $\delta(\mathcal E_T,\mathcal F_T)=0$. Identity on $\mathcal Z$ gives $\delta(\mathcal F_T,\mathcal F_S)\le e$. Apply Proposition~\ref{prop:chain} to the chain $\mathcal E_T,\mathcal F_T,\mathcal F_S,\mathcal E_S$.
\end{proof}
The inequality, sometimes called no-free-transfer, is a consequence of the classical deficiency triangle inequality, not a result outside that theory. The earlier three-term form included $\delta(\mathcal E_T,\mathcal F_T)$, which is identically zero. Thus that term cannot be interpreted as target information loss. Source recovery is what is bounded here.

The term $e$ is per-parameter mismatch. Under a common prior $\pi$, marginal representation mismatch is at most $e$, by convexity of TV. Small marginal mismatch does not imply small $e$; substituting the smaller marginal quantity in the theorem is unjustified. When the target-to-source deficiency is zero, the bound imposes no positive penalty.

\subsection{An exact binary example}
Let $\Theta=\{0,1\}$, let $\mathcal E_S$ observe $\theta$ without error, and let $\mathcal E_T$ observe $\theta$ flipped with probability $a\in[0,1/2]$. Then
\[
 \delta(\mathcal E_S,\mathcal E_T)=0,\qquad
 \delta(\mathcal E_T,\mathcal E_S)=a.
\]
The first equality uses the binary noise channel. For the second, the identity decoder attains error $a$. Any decoder has worst-case error at least its uniform-prior average error, whose minimum for this binary channel is $a$.

Apply to both experiments a channel that flips its input with probability $b\in[0,1/2]$. Source reconstruction error is exactly $b$, by the same argument. The target representation has flip probability $a+b-2ab$, whereas the source representation has flip probability $b$. Hence
\[
 e=a(1-2b),\qquad b+e=a+b(1-2a)\ge a.
\]
For $a=1/4$, the pairs $(b,e)$ at $b=0,1/4,1/2$ are $(0,1/4)$, $(1/4,1/8)$, and $(1/2,0)$. Perfect conditional invariance at $b=1/2$ eliminates all information in either representation. This is a finite illustration with an exact calculation, not an empirical benchmark or evidence of computational hardness.

\section{Constructive Certificates and Optimization Surrogates}\label{sec:optimization}
A constructive simulation claim requires specifying an implementable kernel and bounding its error uniformly over the stated parameter space. For finite experiments, a candidate stochastic matrix can be checked directly. Optimizing statistical deficiency can be formulated as a linear program: introduce stochastic entries $K_{xy}$, absolute-value slack variables $u_{\theta y}$, and minimize $t$ subject to
\[
 \begin{split}
 &K_{xy}\ge0,\quad\sum_yK_{xy}=1,\\
 &u_{\theta y}\ge\sum_xP_\theta(x)K_{xy}-Q_\theta(y),\\
 &u_{\theta y}\ge Q_\theta(y)-\sum_xP_\theta(x)K_{xy},\quad
 \tfrac12\sum_yu_{\theta y}\le t\quad(\theta\in\Theta).
 \end{split}
\]
This computes statistical deficiency for explicit finite tables, up to numerical solver tolerance. It is not a polynomial-time guarantee in the size of an implicitly represented original problem: the tables themselves may be exponentially large. A tabulated solution is also not automatically a uniform efficient simulator across input sizes.

\subsection{A KL certificate under explicit assumptions}
If an implementable $K$ satisfies
$\sup_\theta D_{\rm KL}(P_\theta K\|Q_\theta)\le2\varepsilon^2$,
Pinsker's inequality implies $d(K;\mathcal E,\mathcal F)\le\varepsilon$. Here densities, the KL direction, and the uniform parameter bound must be specified. An empirical average of negative log likelihood is not this KL quantity: it includes an entropy term and sampling error. Nor does an average over parameters imply a supremum bound.

\subsection{Surrogates do not supply the missing bound}
Neural simulator training may use likelihood losses, reconstruction penalties, or maximum mean discrepancy (MMD) \citep{gretton2012kernel}. These can guide optimization, but small surrogate loss alone does not establish small TV or a uniform deficiency certificate. A required additional argument is an explicit inequality connecting the chosen population objective to the desired uniform error, followed by a finite-sample generalization bound and runtime analysis.

For example, on $\{0,1\}$ equip a characteristic kernel with $k(x,y)=c\,\mathbf1\{x=y\}$ for $c>0$. For distinct point masses, TV equals one while MMD equals $\sqrt{2c}$, which can be arbitrarily small as the kernel is rescaled. Any claimed conversion must therefore at least control the kernel and its constants. This example does not rule out conversions under additional fixed-kernel and distributional assumptions.

We make no consistency claim for a generic reconstruction-plus-MMD objective and no assertion that universal approximation alone yields an efficient learned simulator.

\section{Discussion and Limitations}\label{sec:discussion}
The central distinction is between a task diagnostic and a simulation certificate. A certificate bounds a specified kernel's error uniformly over parameters and thereby controls every bounded downstream decision. A diagnostic measured on one task or one sample does not provide that guarantee. For computation, a second distinction separates unconstrained existence from implementation by a uniform machine with a declared resource budget.

These distinctions narrow what can be concluded about approximate solvers. Full-TV approximation to an exact point-mass target demands exact output with high probability. Broader notions of useful approximation require a deliberately chosen coarser target or a restricted decision class. Likewise, a pointwise infimum over polynomial-time machines is not a uniform asymptotic guarantee.

The mathematical bounds in this manuscript are finite-experiment statements. They do not require local asymptotic normality. Estimating their ingredients from data is a separate problem: uniform confidence bounds, simulator selection bias, parameter-space coverage, and computational costs remain to be addressed for particular applications. Passive experiment comparison does not by itself certify interventions; interventional families would have to be specified explicitly.

The exact binary example is an illustration only. We do not claim a validated genomic certificate, neural-solver benchmark, or a new complexity-class separation. The framework's practical value depends on finding targets that express the intended decisions and kernels whose error and implementation cost can both be verified.

\bibliographystyle{plainnat}
\bibliography{references}

\end{document}